\documentclass[11pt]{amsart}
\usepackage{amsmath,amssymb}

\theoremstyle{plain}
\newtheorem{propn}{Proposition}[section]
\newtheorem{thm}[propn]{Theorem}

\theoremstyle{definition}
\newtheorem{defn}[propn]{Definition}

\theoremstyle{remark}
\newtheorem*{rem}{Remark}

\newenvironment{rlist}
{

\begin{enumerate}}
{\end{enumerate}}

\numberwithin{equation}{section}

\begin{document}

\title{Curvature-adapted submanifolds of symmetric spaces}
\author{Thomas Murphy}

\address{School of Mathematical Sciences, University College Cork, Ireland.}

\curraddr{D\'{e}partment de Math\'{e}matique,
Universit\'{e} Libre de Bruxelles,
 \ Boulevard du Triomphe,
B-1050 Bruxelles,
Belgique.}

\email{tmurphy@ulb.ac.be}

\subjclass[2010]{Primary   53C12, 53C40.}

\date{January 3, 2011.}

\keywords{Curvature-adapted hypersurfaces, isoparametric hypersurfaces, complex two-plane Grassmannians, Cayley projective
and hyperbolic planes}

\thispagestyle{empty}

\begin{abstract}
Curvature-adapted submanifolds have been extensively \\ studied in complex and quaternionic space forms. This paper extends
their study to a wider class of ambient spaces. We generalize Cartan's theorem classifying  isoparametric hypersurfaces of spheres  to any compact symmetric space. Our
second objective is to investigate such hypersurfaces in some specific symmetric spaces. We classify those with constant principal curvatures in the Octonionic planes. Various classification results for
hypersurfaces in complex two-plane Grassmannians are also obtained.
\end{abstract}

\maketitle

\section{Introduction}

Let $M$ be a connected hypersurface of a Riemannian manifold $\overline{M}$, $\overline{R}$ be the Riemannian curvature
tensor of $\overline{M}$ and $\xi$ a unit normal vector of $M$ at $p\in M$. The normal Jacobi operator
$$
K_{\xi} := \overline{R}(\xi, \cdot)\xi \in End(T_pM)
$$
of $M$ (with respect to $\xi$) describes the curvature of the ambient manifold $\overline{M}$ at $p$, whereas the shape
operator $A_{\xi}$ of $M$ (with respect to $\xi$) describes the curvature of $M$ as a submanifold of $\overline{M}$ in
direction $\xi$. Both of these are self-adjoint operators, and hence have eigendecompositions. $M$ is said to be curvature
adapted if these operators are simultaneously diagonalizable at every point $p\in M$. This means that a common eigenbasis for $K_{\xi}$ and $A_{\xi}$ exists at every point, which will
generically be denoted by $E$. This condition can be generalized to submanifolds of higher codimension:

\begin{defn}
A submanifold $M$ of $(\overline{M},g)$ is said to be curvature-adapted if the following two conditions are satisfied at
every point $p\in M$:
\begin{rlist}
\item The normal Jacobi operator $\overline{R}(\xi, X)\xi =: K_{\xi}(X)\in T_pM$ for
every unit normal vector field $\xi$ of $M$ and $X\in T_pM$.
\item $A_{\xi}\circ K_{\xi} = K_{\xi}\circ A_{\xi}$, that is the
normal Jacobi and shape operators of $M$ commute.
\end{rlist}
\end{defn}

For a hypersurface $(i)$ is an immediate consequence of the curvature identities. Gray's Theorem \cite{gray} states that any tube around a curvature-adapted submanifold of a locally symmetric space is itself curvature-adapted. Thus  understanding how curvature-adapted hypersurfaces of a given symmetric space arise yields information about the condition in general.

In the real space forms it is easy to see that every submanifold is curvature-adapted.
For other symmetric spaces the condition is restrictive. Curvature-adapted submanifolds in complex and quaternionic space forms has been a particularly fruitful field of study and there is a
substantial body of literature concerned with their classification, of which we just mention \cite{berndt1}, \cite{berndt2}, \cite{datri}, \cite{gray}, \cite{kimura}. Every complex submanifold of a complex space form satisfies this condition, yielding an important family of examples.  In complex space forms the curvature-adapted hypersurfaces coincide precisely with the
Hopf hypersurfaces. These are hypersurfaces with the property that $-J\xi$ is an eigenvector of the shape operator, with
corresponding principal curvature function denoted $\alpha$ (the Hopf principal curvature). Further canonical examples of curvature-adapted submanifolds in general include hyperspheres in any symmetric space and orbits of Hermann actions.

A second family of submanifolds which have been the focus of much attention are isoparametric hypersurfaces.
\begin{defn}
 A smooth function $f: M\rightarrow \mathbb{R}$ is said to \emph{transnormal} if there is a smooth function $b$ such that
$$
\| df\|^2 = b(f).
$$
A transnormal function is said to \emph{isoparametric} if
$$
\Delta f = a(f),
$$
for another continuous function $a$.
\end{defn}

The first equation is equivalent to the level sets $f^{-1}(c)$ being parallel, the second to the level sets having constant
mean curvatures. The level sets are then said to be isoparametric hypersurfaces. Thus if a hypersurface is  isoparametric
all parallel hypersurfaces have constant mean curvature. In space forms these arose naturally in the study of geometrical
optics and their classification here has a long history. Whilst the question is settled for Euclidean and hyperbolic space,
Cartan failed in his attempts to classify them in spheres. Today there is a vast literature on this central problem (see
\cite{thorbergsson} for an excellent survey). Denote by $M_t$ the parallel hypersurfaces at distance $t$ from $M_0=M$, and
the corresponding shape operator at points along the normal geodesic $C_{\xi}(t)$ passing through a point $p\in M_0$ by
$A_{\xi}(t)$. Observe that if the ambient manifold is symmetric the Riccati equation along $C_{\xi}(t)$,
$$
A_{\xi}'(t) = (A_{\xi}(t))^2 + K_{\xi}(t),
$$
simplifies to a family of easy-to-solve ODE's if the hypersurface is curvature-adapted; namely
$$
\lambda_i'(t) = \lambda_i^2(t) + \kappa_i^2
$$
when $\overline{M}$ is compact. Here, and throughout the paper, $\lambda_i$ denotes a principal curvature function and $\kappa_i^2$ an eigenvalue of the normal Jacobi operator with $\kappa_i \geq 0$. There is an analogous formula if $\overline{M}$ is noncompact. This suggests that it is profitable to study this family of hypersurfaces when one wants to analyze the properties of nearby parallel hypersurfaces.

It is a consequence of this equation that if $M\subset G/K$ is a complete curvature-adapted hypersurface with constant principal
curvatures in a rank one symmetric space then it is isoparametric. It is by no means clear if these are the only
isoparametric hypersurfaces. Cartan \cite{cartan} showed that this was indeed the case in the spheres. Wang \cite{wang}
however showed that in general the theory of isoparametric hypersurfaces and that of hypersurfaces with constant principal
curvatures are different in symmetric spaces. Specifically, he discovered families of isoparametric hypersurfaces in
$\mathbb{C}P^n$ with nonconstant principal curvatures. He also showed that a Hopf hypersurface $M\subset \mathbb{C}P^n$ is
isoparametric if and only if it has constant principal curvatures. The obvious question is whether this is true for the
non-compact dual $\mathbb{C}H^n$, or indeed for symmetric spaces generally. This question leads to our first main result:
\begin{thm}\label{thm4}
Let $M$ be a curvature-adapted hypersurface of a compact symmetric space. Then $M$ is isoparametric if, and only if, it has
constant principal curvatures and the eigenvalues of $K_{\xi}$ corresponding to the basis $E$ are constant on $M$.
\end{thm}
In rank one symmetric spaces the eigenvalues of $K_{\xi}$ are constant. Hence for spheres we obtain a new proof of Cartan's theorem.  This
theorem also yields an intrinsic proof of Wang's theorem in the case of $\mathbb{C}P^n$. Wang's proof used the Hopf
fibration $S^{2n+1}\rightarrow \mathbb{C}P^n$ and then applied results from the study of isoparametric functions on spheres.

Our second aim in this paper is to understand the geometry of curvature-adapted hypersurfaces in more general symmetric
spaces than have previously been investigated.  We investigate the geometry of such hypersurfaces in the octonionic projective and hyperbolic planes, denoted $\mathbb{O}P^2$ and $\mathbb{O}H^2$ respectively, where almost nothing is known. Examples analogous to curvature-adapted hypersurfaces in other rank one symmetric spaces will be uncovered and their principal curvatures calculated. We conjecture there are none beyond these examples. Complete curvature-adapted hypersurfaces with
constant principal curvatures are then classified.
\begin{thm}\label{thm2}
Let $M\subset \mathbb{O}P^2$ denote a complete curvature-adapted hypersurface. Then $M$ has
constant principal curvatures if, and only if, $M$ is a principal orbit of a cohomogeneity one action.
\end{thm}
This theorem completes the program initiated in \cite{berndt1}, \cite{berndt2} and classifies all complete curvature-adapted hypersurfaces of compact rank one symmetric spaces with constant principal curvatures. Together with the first Theorem, it also classifies the complete curvature-adapted hypersurfaces of $\mathbb{O}P^2$ which are isoparametric.

Finally we study curvature-adapted submanifolds of complex two-plane Grassmannians
$G_2(\mathbb{C}^{m+2})$. These are rank two Hermitian symmetric spaces with the remarkable property that they are the only
closed non-Ricci-flat Riemannian manifolds which admit both complex and quaternionic-K\"{a}hler structures, denoted $J$ and
$\mathcal{J}$ respectively.

Define $\mathbb{H}\xi(p) := \lbrace J_1\xi(p), J_2\xi(p), J_3\xi(p)\rbrace$, where $J_1,J_2, J_3$ are a local section of
$\mathcal{J}$ at $p\in G_2(\mathbb{C}^{m+2})$. Let $M$ be a real hypersurface of $G_2(\mathbb{C}^{m+2}) =
SU(m+2)/S(U(m)\times U(2))$, and again denote nearby parallel hypersurfaces by $M_t$. There exists an almost Hermitian
structure $J_1\in\mathcal{J}$ and a unit vector $Z\perp \mathbb{H}\xi$ so that at a point $p\in M$
$$
J\xi(p) = \cos(\alpha(p))J_1\xi(p) + \sin(\alpha(p))J_1Z(p),
$$
where $0\leq\alpha(p) \leq\frac{\pi}{2}$. Hypersurfaces such that $\alpha(p)\in \lbrace 0, \frac{\pi}{2}\rbrace$ are
classified in \cite{berndtsuh}: they are precisely the homogeneous hypersurfaces. Moreover, it can be calculated that  they
are curvature-adapted. We conjecture there are no more curvature-adapted hypersurfaces in $G_2(\mathbb{C}^{m+2})$. As
evidence of this, we show that there generically are no curvature-adapted hypersurfaces satisfying one additional
constraint.

\begin{thm}\label{thm1}
There are no curvature-adapted hypersurfaces of $G_2(\mathbb{C}^{m+2})$  such that  $\cos(\alpha(p))$ $\notin \lbrace 0,\frac{3}{5}, \frac{4}{5}, 1 \rbrace$ $\forall p\in M$ and either $\langle A_{\xi}(J_1\xi(t)), J_1\xi(t)\rangle$,
or $\langle A_{\xi}(J_1Z(t)), J_1Z(t)\rangle$, or their ratio is constant along the normal geodesic $C_{\xi}(t)$ through all
points $p\in M$.
\end{thm}

In cases where $\cos(\alpha) \in \lbrace \frac{3}{5}, \frac{4}{5}\rbrace $ the eigenspaces of $K_{\xi}$ change and have
different multiplicity, and our approach fails. However, we hope this result with stimulate further work on the study of curvature-adapted submanifolds in higher rank symmetric spaces.

\section{Isoparametric hypersurfaces}
Historically, there have been two  definitions of isoparametric hypersurfaces of $\overline{M}$. Cartan \cite{cartan}, following Somilgiania,
Levi-Civita and Segre, defined a hypersurface $M\subset \overline{M}$ to be isoparametric if $M\simeq f^{-1}(t)$, where $f:
\overline{M}\rightarrow \mathbb{R}$ is an isoparametric function. We will follow this definition. For the alternative, which we define as \emph{weakly isoparametric},  a hypersurface $M\subset \overline{M}$ is required to have all constant mean curvatures for all sufficiently close parallel hypersurfaces $M_t$, $|t| < \epsilon$.  For spheres and complex projective spaces it is known that
any weakly isoparametric hypersurface is an open part of the level set of an isoparametric function, so there the
definitions are equivalent. We also remark that the exceptional orbits of the Riemannian foliation induced by an
isoparametric function on any Riemannian manifold are minimal submanifolds (\cite{getang}, \cite{wang}): this is one of the
few general methods of constructing minimal submanifolds of a Riemannian manifold. Throughout this paper $M_{\xi}$
refers to the focal set of a hypersurface $M\subset \overline{M}$, and connected components of the focal set are denoted
$Q_i$.

Curvature-adapted hypersurfaces were introduced by d'Atri \cite{datri}, who observed that all known examples of
isoparametric hypersurfaces in rank-one symmetric spaces are curvature-adapted. He also generalized Cartan's fundamental
formula for isoparametric hypersurfaces in spheres to any rank one symmetric space. We remark that an elegant proof of his theorem may be deduced from the following two steps: firstly using that $M$ is isoparametric to deduce that the focal manifolds are minimal \cite{wang2}, and secondly
calculating the principal curvatures of $M$ in terms of the principal curvatures of $Q_1$.

The Riccati equation allows much to be said for curvature-adapted hypersurfaces and the relationships between their geometry
and that of nearby parallel hypersurfaces. As such, one expects curvature-adapted hypersurfaces and isoparametric
hypersurfaces in symmetric spaces to be intimately related. Our next objective is to justify this statement by establishing
Theorem \ref{thm4}.

\proof Let $M^n\subset G/K$ be a complete curvature-adapted hypersurface. If it has constant principal curvatures and the eigenvales
of $K_{\xi}$ associated to $E$ are constant, it is an easy consequence of the Riccati equation that $M$ is isoparametric.

Conversely, suppose that $M$ is isoparametric, but assume that either: (i) the eigenvalues of $K_{\xi}$ with respect to $E$
are nonconstant, (ii) the principal curvatures of $M$ are nonconstant, or (iii) both the principal curvatures of $M$ and the
eigenvalues of $K_{\xi}$, denoted $\kappa_i$ with respect to $E$ are nonconstant. Only the proof of (iii)  will be
given; all other cases are analogous. Suppose that the principal curvatures of $M$ are nonconstant. Let $p\in M$ be a fixed
point. Assume there is a point $q\in M$ where the principal curvatures differ from $p$. The strategy of the proof is to
consider the Riccati equation at these two points. Since $M =: M_0$ is isoparametric, the sum $\sum_i\lambda_i(0) = c_0$,
and similarly for all parallel hypersurfaces $M_t$, $|t|<\epsilon$, one has $\sum_i \lambda_i(t) = c_t$. Here $\epsilon$ is
$sup\lbrace |t|: M_t \text{ is a hypersurface} \rbrace$, and for ease of notation we will assume that $\lambda_i(p)(0) \neq
\lambda_j(p)(0)$ for all $i\neq j$. If any principal curvature has multiplicity greater than one, the same proof goes
through with some trivial modifications. Therefore $\sum_i\lambda_i(p)(t) = \sum_i\lambda_i(q)(t)$ for all $t, |t|<
\epsilon$. But from the Riccati equation one may solve to obtain $$\lambda_i(p)(t) = \kappa_i\cot(\theta_i(p) -
\kappa_i(p)t),$$ where $\kappa_i\cot(\theta_i) = \lambda_i$, and similarly at $q$. Hence
$$
\sum_{i=1}^n \kappa_i(p)\cot(\theta_i(p) -\kappa_i(p)t) - \sum_{i=1}^n\kappa_i(q)\cot(\theta_i(q) -\kappa_i(q)t) = 0
$$
for $|t|< \epsilon$.  Expanding out the Taylor expansion for $\cot$ around $ t=0$ and grouping coefficients yields a
polynomial $F(t) =0$ for all $|t|<\epsilon$. This cannot vanish unless all coefficients of $F$ vanish. This is equivalent
to
$$
\sum_{i=1}^n \frac{1}{r_i(p) -t} - \sum_{i=1}^n\frac{1}{r_i(q) - t}- \sum_{j = 0}^{\infty} f_jt^j = 0,
$$
where $r_{i} = \frac{\theta_i}{\kappa_i}$ and $f_j\in C^{\infty}(\mathbb{R})$, for all $t, |t|<\epsilon$.

Let $|r_1(p)| = min \lbrace |r_i(p)|, |r_i(q)|\rbrace$.  Suppose there is a focal set, so $|r_1(p)| = \epsilon$. Multiplying
across by $\left( r_1(p)-t\right)$ and taking the limit as $t\rightarrow r_1(p)$ yields a contradiction (since $\sum_{j =
0}^{\infty} f_jt^j$ is convergent in a compact neighbourhood of $r_1(p)$) unless $r_1(p)=r_1(q)$. Repeating this argument if
necessary shows that for all the $k$ principal curvatures which focalize the corresponding term $r_i$ is constant. Thus the
isoparametric condition may be rewritten as
\begin{equation*}
\frac{1}{r_k(p)-t} = \sum_{i=k+1}^n \frac{1}{r_i(p) -t} - \sum_{i=k}^n\frac{1}{r_i(q) - t}- \sum_{j = 0}^{\infty}
\tilde{f}_jt^j,
\end{equation*}
where $r_{i} = \frac{\theta_i}{\kappa_i}$, for all $t, |t|<\epsilon$.  Then suppose  $|r_k(p)|$ $= min \lbrace |r_i(p)|,
|r_i(q)|, i=k,\dots, n\rbrace$. Since both sides of this equation are infinite power series with constant coefficients, if
they agree on an open interval (namely $|t|<\epsilon$) then all coefficients must agree and so they agree on any interval
for which both power series converge. This implies that we can repeat the above argument and take the limit as $t\rightarrow
r_2(p)$ to again derive a contradiction unless $r_2(p)=r_2(q)$, and so forth.

Repeating this argument yields $r_i(p)=r_i(q)$ for all $i$. Suppose without loss of generality $\kappa_i(p) > \kappa_i(q)$, so that
$\theta_i(p) < \theta_i(q)$. Notice $\theta_i(p)$ and $\theta_i(q)$ have the same sign; we assume without loss of generality they are both
negative. Then considering nearby hypersurfaces $M_t$ with $t>0$ we see from the solution to the Riccati equation that for some $t_0>0$ one has
$\tilde{\theta}_i(p)(t_0) > 0 > \tilde{\theta}_i(q)(t_0)$, where $$ \tilde{\theta}_i(p)(t_0) = (\theta_i(p) - \kappa_i(p)t_0) $$ and similarly
at $q$. Repeating this calculation with $M = M_{t_0}$ if necessary gives a contradiction, and the proof is complete. \endproof

\begin{rem}
The same proof yields an analogous result  in the non-compact case with $\coth$ replacing $\cot$, under the assumption that
$|\lambda_i(p)|\geq |\kappa_i(p)|$ $\forall i$.
\end{rem}

Suppose there exists a real number $\epsilon>0$ such that all hypersurfaces $M_t$ within distance $ |t| < \epsilon$ of a
given hypersurface $M$ of a symmetric space have constant mean curvatures. Such hypersurfaces are said to be \emph{ weakly
isoparametric}.   Let $\mathbb{K}^n(c), c\neq 0$, $\mathbb{K} = \mathbb{C}$ or $\mathbb{H}$ denote the simply connected manifolds with constant holomorphic (resp. quaternionic) sectional curvature $c$. For such spaces we obtain a
stronger result than above: the classification of weakly isoparametric curvature-adapted hypersurfaces.
\begin{thm}\label{j}
Let $M\subset \mathbb{K}^n(c), c\neq 0$, $\mathbb{K} = \mathbb{C}$ or $\mathbb{H}$ be a curvature-adapted hypersurface. Then $M$ is weakly isoparametric if and only if it has
constant principal curvatures.
\end{thm}
The proof is an adaptation of a proof due to Cecil-Ryan \cite{cecilryan} for isoparametric hypersurfaces in spheres. We note
that for $\mathbb{C}P^n$ this result was already proved by Wang using the Hopf fibration, but again our approach has the
advantage of giving an intrinsic proof.

\proof We just give the proof for $\mathbb{K}= \mathbb{C}$, the case of $\mathbb{H}$ is analogous. Let $M\subset \mathbb{C}^n(c), c\neq 0$ be a Hopf hypersurface. It is obvious from the work of Berndt \cite{berndt}-
Kimura\cite{kimura} that if $M$ has constant principal curvatures it is weakly isoparametric. Suppose conversely $M$ is
weakly isoparametric but does not have constant principal curvatures. Since $M$ is weakly isoparametric, the sum of
principal curvatures $\sum_{i=1}^{2n-1}\lambda_i(t)$ is constant on $M_t$. It is well known that the principal curvature
$\lambda_1(t) = \alpha(t)$ corresponding to $-J\xi(t)$ is constant (\cite{ki},\cite{maeda}): subtracting this from the
equation yields that $\sum_{i=2}^{2n-1}\lambda_i(t)$ is constant on $M$(t). Now, observing that all the terms in this
equation correspond to an eigenvector of the normal Jacobi operator with eigenvalue $\pm1$, we can adopt the proof given in
\cite{cecilryan}. Differentiating this equation with respect to $t$ yields that $\sum_{i=2}^{2n-1}\lambda_i^2(t)$ is
constant on $M_t$. Differentiating this again forces
$$
\frac{-c}{2}\sum_{i=2}^{2n-1}\lambda_i(t) + 2\sum_{i=2}^{2n-1}\lambda_i^3(t)
$$
and hence $\sum_{i=2}^{2n-1}\lambda_i^3(t)$ to be constant on  $M_t$. Iterating this calculation $n$ times, one obtains that
$\sum_{i=2}^{2n-1}\lambda_i^k(t)$ is constant on $M_t$ for $k=1,\dots,n$. But it is known $\alpha = \lambda_1$ is constant.
From Newton's identities it follows that $\lambda_i(t)$ are constant on $M_t$, $i=1,\dots, 2n-1$, and we are done.
\endproof

\section{The octonionic projective and hyperbolic planes}

The octonionic (or  Cayley)  projective plane $\mathbb{O}P^2$ and its noncompact dual $\mathbb{O}H^2$ are intriguing
mathematical objects. As rank one symmetric spaces one would expect their geometries to be well-understood, yet they remain
mysterious objects. Due to their esoteric nature a brief exposition of their properties is presented. We refer the
interested reader to \cite{baez} for an in-depth study. $\mathbb{O}$ will denote the octonions, an eight dimensional
non-associative division algebra over $\mathbb{R}$ which satisfies the alternative law. $\mathbb{O}$ has a multiplicative
identity $1$ and a positive definite bilinear form $\langle, \rangle$ whose associated norm $||, ||$ satisfies $ ||ab|| = ||a||.||b||$. As is the case for the complex numbers and quaternions, each element $a\in \mathbb{O}$ can be expressed in the
form $ a = \alpha1 + a_0$ where $\alpha\in\mathbb{R}$ and $\langle \alpha, a_0 \rangle = 0$. A conjugation map is defined as
$a\rightarrow a^* = \alpha1 - a_0$. This is an anti-automorphism: $(ab)^* = b^*a^*$. To measure nonassociativity there is
the associator
$$
(a,b,c) := (ab)c - a(bc).
$$
Linearization of the associative law $(a,a,c) = (c,a,a) = 0$ yields that $(a^*,b, c) = -(a,b,c)$. This implies that $
(a,a^*,b) = -(a^*,a^*,b) =0$. We define a canonical basis of $\mathbb{O}$ to be any basis of the form $\lbrace 1 = J_0, J_1,
\dots, J_7\rbrace$ such that
\begin{rlist}
\item $\langle J_i,J_j\rangle = \delta_i^j$,
\item $J_i^2 = -1, i\neq 0$,
\item $J_iJ_j + J_jJ_i = 0, i\neq j, i,j\neq 0$,
\item $J_iJ_{i+1} = J_{i+3}$ ,modulo $7, i\neq 0.$
\end{rlist}

The Cayley projective and hyperbolic planes have an \emph{octonionic} structure pointwise. We will outline this construction
for $\mathbb{O}P^2$. Firstly $\mathbb{O}P^2= F_4/Spin(9)$ is a sixteen dimensional rank one symmetric space. Therefore the isotropy
representation of $Spin(9)$ acts irreducibly on $T_p\mathbb{O}P^2 = \mathbb{R}^{16}$. This induces a transitive group action
on the sphere $S^{15}\subset T_p\mathbb{O}P^2$. But $S^{15}=Spin(9)/Spin(7)$. So $Spin(7)$ fixes a point on $S^{15}$, or
equivalently some vector $X\in T_p\mathbb{O}P^2$. If we restrict the $Spin(9)$ action to $Spin(7) \subset Spin(9)$, then
$Spin(7)$ fixes $X$, and so acts trivially on $\mathbb{R}X\subset T_p\mathbb{O}P^2$. It therefore leaves invariant two
subspaces of $\mathbb{R}^{15}\bot \mathbb{R}X$, namely $\mathbb{R}^7$ and $\mathbb{R}^8$. The induced representations of
$Spin(7)$ are the standard representation and the spin representation respectively. We now identify $T_p\mathbb{O}P^2$ with
$\mathbb{O}\oplus\mathbb{O}$ by equating one copy of the octonions with $\mathbb{R}\oplus\mathbb{R}^7$, and the second copy
of $\mathbb{O}$ with $\mathbb{R}^8$.

Both $\mathbb{O}P^2$ and $\mathbb{O}H^2$ have unique $Ad(Spin(9))$-invariant Riemannian metrics up to homothety. The
curvature tensors for these manifolds are very different to the other compact rank one symmetric spaces. It was written down
for the first time in \cite{browngray} for $\mathbb{O}P^2$ (resp $\mathbb{O}H^2$) at
$T_p\overline{M}=\mathbb{O}\oplus\mathbb{O}$ as
\begin{align*}\label{tm4}
R((a,b),(c,d))(e,f) =& \frac{\pm1}{4}(\lbrace 4\langle c,e \rangle a -4\langle a,e \rangle c +(ed)b^*\\ &- (eb)d^* +
(ad-cb)f^*\rbrace, \lbrace 4\langle d,f \rangle b - 4\langle b,f \rangle d\\ &+ a^*(cf) - c^*(af) - e^*(ad-bc)\rbrace).
\end{align*}
The inner product $\langle,\rangle_p$ induced by the metric is given by
$$
\langle (a,b), (c,d)\rangle = \langle a,c\rangle + \langle b,d\rangle.
$$
Throughout, we scale the metric to have sectional curvatures to lie between $\pm1$ and $\pm4$.

Let $M$ be a real hypersurface of $\mathbb{O}P^2$ with normal vector field $\xi$. Then along the normal geodesic $C_{\xi}$
we can parallel translate our basis of $T_pM$, $\mathbb{O}\oplus\mathbb{O}$, and hence get an invariant description of the
Riemannian curvature tensor along $C_{\xi}$ in terms of $\mathbb{O}\oplus\mathbb{O}$. This is because the isotropy group of
the geodesic $C_{\xi}$ is $Spin(7)\subset Spin(9)$, and as we have seen the induced action of $Spin(7)$ on
$T_p\mathbb{O}P^2= \mathbb{R}^8 \oplus\mathbb{R}^8$ decomposes into actions on $\lbrace \xi \rbrace \oplus
\mathbb{R}^7\oplus\mathbb{R}^8$. The $\mathbb{R}^7$ is calculated to correspond to the $+4$ eigenspace of $K_{\xi}$ and the
second copy of $\mathbb{R}^8$ corresponds to the $+1$ eigenspace of $K_{\xi}$. Taking $\xi = (1,0)$ along $C_{\xi}$,  we may
choose as our common eigenframe along the geodesic $C_{\xi}$ the basis
$$
E(t) = \lbrace U_i(t), V(t), J_iV(t)\rbrace,
$$
where $J_i$ is the octonionic structure along $C_{\xi}$, $U_i = -J_i\xi$ span the $\pm1$ eigenspace of $K_{\xi}$ and
$V,J_iV$ span the $\pm 4$ eigenspace of $K_{\xi}$. It is to be assumed we are working with this frame in what follows.

Such structures only exist along the geodesic $C_{\xi}$: there cannot be a parallel rank seven subbundle $\mathcal{J}\subset
End(T\overline{M})$, $\overline{M} = \mathbb{O}P^2$ or $\mathbb{O}H^2$. This can be seen by the following argument, shown to
us by Robert Bryant. If such a bundle were to exist, its holonomy would have to be a quotient group of $Spin(9)$ that can be
embedded in $SO(7)$, but the only such subgroup (since the Lie algebra of $Spin(9)$ is simple and $Spin(9)$ is connected) is
the trivial subgroup. Thus, if there were such a bundle, it would have a basis of Levi-Civita parallel sections.  In
particular, the action of $Spin(9)$ on $T_p\overline{M}$ would have to commute with all of these endomorphisms, meaning that
the space of endomorphisms of $T_p\overline{M}= \mathbb{R}^{16}$ which commute with the $Spin(9)$ action would have
dimension at least $7$. However the space of linear transformations which commute with $Spin(9)$ on $\mathbb{R}^{16}$ is
one-dimensional, so we derive a contradiction. To see this, suppose that a linear transformation $T$ of $\mathbb{R}^{16}$
commutes with $Spin(9)$. Thus $T$ acts on $\mathbb{R}^{16}$ commuting with the Lie algebra $spin(9)$, and so on
$\mathbb{C}^{16}$ commuting with the Lie algebra $spin(9,\mathbb{C})=spin(9)+Ispin(9)$, where $I$ denotes the complex
structure. Each eigenspace of $T$ must be $Spin(9,\mathbb{C})$-invariant. But $Spin(9,\mathbb{C})$ acts irreducibly on
$\mathbb{C}^{16}$, so $T$ must have the whole of $\mathbb{C}^{16}$ as eigenspace. In other words, $T$ acts as rescaling by a
single complex number. $T$ is real, so its eigenvalues come in complex conjugate pairs and so there is a single real
eigenvalue. Hence $T$ is a real rescaling.

\begin{propn}\label{bob}
The following are curvature-adapted hypersurfaces;
\begin{rlist}
\item  the tube of radius $r$ around a
totally geodesic $\mathbb{O}P^k\subset \mathbb{O}P^2$, where $r\leq \frac{\pi}{4}$ and $0\leq k <2$,\\
\item the tube of radius $r$ around a totally
geodesic $\mathbb{H}P^2\subset \mathbb{O}P^2$, where $0<r\leq \frac{\pi}{4}$,\\
\item the tube of radius $r\in\mathbb{R}^+$ around a totally geodesic
$\mathbb{O}H^k\subset \mathbb{O}H^2$, where $0\leq k <2$,\\
\item the tube of radius $r\in\mathbb{R}^+$ around a totally
geodesic $\mathbb{H}H^2\subset \mathbb{O}H^2,$\\
\item a horosphere in $\mathbb{O}H^2$.\\
\end{rlist}
Their principal curvatures, together with their multiplicities, are given  in the following table;

\bigskip

\begin{tabular}
{| l | c | c | c | c| r | } \hline
M &(i) & (ii)& (iii) & (iv) & (v) \\ \hline
$\lambda_1$ & - & $\cot(r)$  & - & $\coth(r)$& $1$ \\
$\lambda_2$ & $-\tan(r)$& $-\tan(r)$ & $\tanh(r)$ & $\tanh(r)$  & - \\
$\alpha_1$ & $2\cot(2r)$ & $2\cot(2r)$ & $2\coth(2r)$ & $2\coth(2r)$ & $2$ \\
$\alpha_2$ & - & $-2\tan(2r)$ & - & $2\tanh(2r)$ & - \\
$m(\lambda_1)$ & - & $4$ & - & $4$ & $8$\\
$m(\lambda_2)$ & $8$ & $4$ & $8$ & $4$ &-\\
$m(\alpha_1)$ & $7$ & $3$ & $7$ & $3$ & $7$\\
$m(\alpha_2)$ & - &$4$ & - &$4$ & - \\
\hline
\end{tabular}
\end{propn}

\proof The proof will firstly be outlined in the first two cases: the noncompact cases are exactly analogous. The proof for
the horosphere follows the same idea of the proof of Theorem 2 in \cite{berndt}, so we omit it.  Fix a point $p$ in
$\mathbb{O}P^1$ and a unit normal vector $\xi$ at p. The $4$-eigenspace of $K_{\xi}$ is seven-dimensional and equal to the
orthogonal complement of $\xi$ in the normal space of $\mathbb{O}P^1$ at $p$.  The $1$-eigenspace is eight-dimensional and
is equal to the tangent space of $\mathbb{O}P^1$ at p. This tells us that $\mathbb{O}P^1 \subset \mathbb{O}P^2$ is curvature
adapted. It follows that the tubes around $\mathbb{O}P^1$ are curvature-adapted. These tubes are the principal orbits of the
action of $Spin(9)$ on $\mathbb{O}P^2 = F_4/Spin(9)$. This is a cohomogeneity one action with two singular orbits; a totally
geodesic $\mathbb{O}P^1 = S^8$ and a single point. Alternatively, choosing this point it is curvature-adapted, and so the
tubes around this point (geodesic hyperspheres) are curvature-adapted. Hence the tubes around $\mathbb{O}P^1\subset
\mathbb{O}P^2$ are nothing more than geodesic hyperspheres.

For $\mathbb{H}P^2\subset \mathbb{O}P^2$, note that $\mathbb{H}P^2$ is an orbit of the maximal subgroup $Sp(3)Sp(1)$ in
$F_4$. However $Sp(1)$ centralizes $Sp(3)$, and so we can restrict to $Sp(3)$. The stabilizer at a point in $\mathbb{H}P^2$
is $Sp(2)Sp(1)$, and the slice representation is the standard representation of $Sp(2)Sp(1)$ on $\mathbb{H}^2 =
\mathbb{R}^8$. Fixing a normal vector $\xi$, the corresponding stabilizer of $Sp(2)Sp(1)$ is $Sp(1)Sp(1)Sp(1)$, and the
restriction of the slice representation to this subgroup decomposes into $\mathbb{H}\oplus\mathbb{H}$. We have seen that the
eigendecomposition of $K_{\xi}$ is $ \lbrace \xi \rbrace \oplus \lbrace U_i \rbrace\oplus \lbrace J_iV \rbrace $. The
intersection with the normal space of $\mathbb{H}P^2$ is
$$
\mathbb{H}\oplus\mathbb{H} = \mathbb{R} \oplus \mathbb{R}^3 \oplus \mathbb{R}^4.
$$
This tells us that $K_{\xi}:T\mathbb{H}P^2\rightarrow T\mathbb{H}P^2$, and it follows $\mathbb{H}P^2$ is curvature-adapted,
as is the tube around it.

We will calculate the principal curvatures for the tube around $\mathbb{H}P^2$. Analogous techniques yields the principal
curvatures of all the other examples given. Choose a point $q$ lying at distance $r$ from $p$ along the geodesic $C_{\xi}$
determined by $\xi$. Apply again the standard theory to calculate the principal curvatures. Since $\mathbb{O}P^2$ is a
symmetric space, the curvature tensor is parallel, so the eigenvalues of $K_{\xi}$ are constant. From the theory of matrix
differential equations this is the same as solving a family of second order ODE's, the first two of which (those involving
$X_1,X_2\in\mathbb{H}\times 1$ tangent to $\mathbb{H}P^2$ are $ Y_i'' + 4Y = 0, $ with $i=1,2$ and initial conditions
$Y_1(0)=1, Y_1'(0)=0$. The solution is $\cos(2r)$, and the corresponding eigenvalue $-2\cot(2r)$. Then the second two
equations (involving $X_3,X_4\in\mathbb{H}\times 1$ normal to $\mathbb{H}P^2$  arise from solving $ Y_i'' + 4Y = 0, $ with
$i=3,4$ and initial conditions $Y_1(0)=0, Y_1'(0)=1$. Solving yields $Y(r) = \sin(2r)$ and the corresponding principal
curvature is $-2\tan(2r)$. All principal curvatures are calculated using the same method.
\endproof

We conjecture that these are all the curvature-adapted hypersurfaces in $\overline{M}$, which is a completely analogous
situation to the classification in the quaternionic space forms. Establishing this conjecture would complete the
classification of such hypersurfaces in rank one symmetric spaces. To provide some evidence for this conjecture,  let us give the proof of Theorem \ref{thm2}.

\proof Let $M$ be a complete curvature-adapted hypersurface of $\mathbb{O}P^2$ or $\mathbb{O}H^2$ with constant principal curvatures. Then it is not hard to see that $M$ must arise as a level set of an isoparametric function. In this
situation the focal set has two disjoint connected components
$$M_{\xi} = Q_1 \cup Q_2.
$$
For convenience scale the ambient metric to have ambient sectional curvature lying between $1$ and $ 4$ respectively. $M$
may be viewed as a tube around one of its focal manifolds, $Q_1$. As $M$ is curvature-adapted, so is $M_t$ for all $t$.
Hence we may choose $U_i\in E(t)$ for all $t$ and by continuity it follows that at the focal manifold $Q_1=M_{t_0}$ we have
$U_i\in E(t_0)$, where $E(t)$ denotes a common eigenbasis of $A_{\xi}(t)\oplus Id$ and $K_{\xi}$. By similar reasoning for
every $V\in E(t_0)$, one may choose $J_iV\in E(t_0), i=1,\dots, 7$ to complete the common eigenbasis that we will work with
for this proof. The idea of the proof is to show that either $Q_1$ or $Q_2$ is totally geodesic. The tubes around totally
geodesic submanifolds which are curvature-adapted are precisely the homogeneous hypersurfaces.

$M$ is equifocal, so a theorem of Tang \cite{tang} implies $g=1,2$, where $\pi/2g$ denotes the length of the interval
between $Q_1$ and $Q_2$. $g=1$ is immediately ruled out, as any non-zero principal curvature functions on $Q_1$ would
focalize before $t=\frac{\pi}{2}$ unless $Q_1$ or $Q_2$ is totally geodesic. Hence it may be assumed that the distance between $Q_1$ and $Q_2$ is $\frac{\pi}{4}$.
Write $Cot(\theta_i)$ for the nonezero principal curvatures of $Q_1$ at $p$ with respect to the normal vector $\xi(p)$. The solutions of the Riccati equation fall into two families, for the $+4$-eigenspace of $K_{\xi}$ one obtains $ \lambda_i(t) = 2\cot(\theta_i - 2t)$ and for the $+1$-eigenspace $\cot(\theta_i - t)$ with $0\leq t \leq \frac{\pi}{4}$. Corresponding to $\lambda_i(0) = 0$ principal curvatures at $Q_1$ the corresponding solutions are $-2\tan(2t)$ and $-\tan(t)$ respectively.

 From this it follows all the principal curvatures of $Q_1$ and $Q_2$ in the $+4$-eigenspace must be zero. Moreover, the only possiblities for the principal curvatures of $Q_1$ in the $+1$-eigenspace are $1, 0$ or $-1$. To see this, observe that each $\theta_i$ is either $\frac{\pi}{4}, \frac{\pi}{2}$, or $\frac{3\pi}{4}$ because by assumption the distance between $Q_1$ and $Q_2$ is $\frac{\pi}{4}$. This can be seen by travelling along the geodesic $C_{\xi}$ which passes through $p$ again at distance $\pi$. Suppose $Q_1$ is not totally geodesic. Then $Q_2$ is not minimal if it is not totally geodesic, a contradiction. \endproof

\section{Complex two-plane Grassmannians}
Consider the $(m+1)$-dimensional $\mathbb{C}P^{m+1}$ embedded canonically as a totally geodesic submanifold of
$\mathbb{H}P^{m+1}$. The focal set $Q^{m+1}$ of $\mathbb{C}P^{m+1}$ is a submanifold of codimension three. At each point of
$Q^{m+1}$ the null space of the shape operator $A_{\xi}$ is independent of the choice of normal vector field $\xi$. It
determines a one-dimensional Riemannian foliation $\mathcal{F}$ on $Q^{m+1}$ by closed geodesics (in both the ambient metric
of $\mathbb{H}P^{m+1}$ and the induced metric on $Q^{m+1}$). The orbit space $B^{m+1} := Q^{m+1}/\mathcal{F}$, equipped with
the Riemannian structure for which the canonical projection $\pi_B$ is a Riemannian submersion is isometric to the
Riemannian symmetric space $(G_2(\mathbb{C}P^{m+2}), \langle, \rangle)$. This fibration yields all the geometric information
about $G_2(\mathbb{C}P^{m+2})$ in terms of the intrinsic and extrinsic structure of the focal set $Q^{m+1}$ of
$\mathbb{C}P^{m+1}$ in $\mathbb{H}P^{m+1}$. $G_2(\mathbb{C}P^{m+2})$ has both a K\"{a}hler structure, $J$ and a quaternionic
K\"{a}hler structure $\mathcal{J}$ induced from the ambient quaternionic K\"{a}hler structure. Let $J_1, J_2,J_3$ denote a
canonical local basis of $\mathcal{J}$ near $p\in M$. Recall that we defined $\mathbb{H}\xi = Span\lbrace \xi, J_i\xi_p :
J_i\in \mathcal{J}_p \rbrace$. Similarly define $\mathbb{C}\xi = Span\lbrace \xi, J\xi \rbrace$. Then the Riemannian
curvature tensor is given as
\begin{align*}
R(X,Y,Z) = &\langle Y,Z\rangle X - \langle X,Z \rangle Y\\
& + \langle JY,Z \rangle JX - \langle JX,Z \rangle JZ - 2 \langle JX,Y \rangle JZ\\
& + \sum_{\nu = 1}^3 \langle J_{\nu}Y,Z \rangle J_{\nu}X - \langle J_{\nu}X,Z \rangle J_{\nu}Z - 2 \langle J_{\nu}X,Y
\rangle J_{\nu}Z\\
& + \sum_{\nu = 1}^3 \langle J_{\nu}JY,Z \rangle J_{\nu}JX - \langle J_{\nu}JX,Z \rangle J_{\nu}JY
\end{align*}

It follows from the expression for the Riemannian curvature tensor that all homogeneous hypersurfaces in
$G_2(\mathbb{C}^{m+2})$ are curvature-adapted. These hypersurfaces are
\begin{rlist}
  \item a principal orbit of the action of $S(U(m+1)\times
  U(1))\subset SU(m+2)$, or
  \item a principal orbit of the action of $Sp(n+1)\subset SU(2n+2)$
  if $m=2n$.
\end{rlist}
In both cases it is possible to calculate that these hypersurfaces are curvature-adapted. For the principal orbits of the
action of $S(U(n+1)\times U(1))$, consider the totally geodesic $\mathbb{G}_2(\mathbb{C}^{n+1})$ arising as an exceptional
orbit of this action. This is known to be complex with respect to the K\"{a}hler structure and quaternionic with respect to
the quaternionic-K\"{a}hler structure \cite{klein}. Hence, a short calculation using the explicit expression for the
Riemannian curvature tensor yields that $K_{\xi}(T\mathbb{G}_2(\mathbb{C}^{n+1})) \subset \mathbb{G}_2(\mathbb{C}^{n+1})$,
whence $\mathbb{G}_2(\mathbb{C}^{n+1})$ and the tubes around it are curvature-adapted. The case where $M$ is a principal
orbit of the action of $Sp(n+1), n=2m$ is completely analogous.

Note that both of these families of hypersurfaces may be viewed as tubes of a fixed radius around maximal totally  geodesic
submanifolds: the second family may be viewed as tubes around a totally geodesic $\mathbb{H}P^n$. Such tubes have at most
five distinct principal curvatures, all of which are constant. Therefore the known curvature-adapted hypersurfaces in
$G_2(\mathbb{C}^{m+2})$ exhibit similar behaviour as in quaternionic projective and hyperbolic spaces. We now prove Theorem
\ref{thm1}:

\proof  We will assume without loss of generality that $\langle A_{\xi}(J_1\xi(t)), J_1\xi(t)\rangle$ is constant along
$C_{\xi}$: under any of the other assumptions given the proof is analogous. The eigenvalues of $K_{\xi}$ together with their
eigenspaces and dimensions fall into one of three possibilities, all listed in tables in \cite{berndt2}. Suppose we are in
the third case of \cite{berndt2}, where there exists an almost Hermitian structure $J_1\in\mathcal{J}$ and a unit vector
$Z\perp \mathbb{H}\xi$ so that $ J\xi = \cos(\alpha)J_1\xi + \sin(\alpha)J_1Z$, where $0<\alpha <\frac{\pi}{2}$ and moreover
suppose that $\cos(\alpha) \notin \lbrace \frac{3}{5}, \frac{4}{5}\rbrace$.  In what follows we will abuse notation slightly
by dropping references to the point $p$. Setting $\beta = \frac{\alpha}{2}$, the eigenspaces corresponding to the
eigenevectors $ X_1= \mathbb{R}\cos(\beta)J_1\xi + \sin(\beta)J_1Z $ and $ X_2=\mathbb{R}\sin(\beta)J_1\xi - \cos(\beta)J_1Z$
are seen to both have dimension one. The corresponding eigenvalues are $-4(1 + \cos(\alpha))$ and  $-4(1 - \cos(\alpha))$.  As
there is a common eigenbasis $E$ for $A_{\xi}$ and $K_{\xi}$ at each point $p$, these vectors must both be eigenvectors of
$A_{\xi}$. Let $A_{\xi}(X_1) = \lambda_1(X_1)$, etc.  Consider the equation
\begin{equation}\label{tm}
\langle A_{\xi} (J_1\xi),
 J_1Z\rangle = \langle J_1\xi, A_{\xi} (J_1Z)\rangle.
\end{equation}
 Rewriting the left hand side gives
\begin{align*}
\langle A_{\xi}(J_1\xi), J_1Z\rangle & = \langle A_{\xi} \Big(\frac{-\sin(\beta)}{\cos(\beta)}J_1Z +
\frac{1}{\cos(\beta)}X\Big), J_1Z\rangle \\
&= \frac{-\sin(\beta)}{\cos(\beta)}\langle A(J_1Z),J_1Z \rangle + \frac{\lambda_1\sin(\beta)}{\cos(\beta)}
\end{align*}
Rewriting the right hand side of the equation produces
\begin{align*}
\langle J_1\xi, A_{\xi}(J_1Z)\rangle &= \langle J_1\xi, A_{\xi}\Big(\frac{1}{\sin(\beta)}X -
\frac{\cos(\beta)}{\sin(\beta)}J_1\xi\Big)\rangle \\
&= \frac{\lambda_1\cos(\beta)}{\sin(\beta)} - \frac{\cos(\beta)}{\sin(\beta)}\langle A(J_1\xi),J_1\xi \rangle.
\end{align*}
Rearranging for $\lambda_1$ yields that
\begin{align*}
\lambda_1\Big(\frac{\sin(\beta)}{\cos(\beta)}-\frac{\cos(\beta)}{\sin(\beta)}\Big)
= &\frac{\sin(\beta)}{\cos(\beta)}\langle A_{\xi}(J_1Z),J_1Z \rangle\\
& -\frac{\cos(\beta)}{\sin(\beta)}\langle A_{\xi}(J_1\xi), J_1\xi\rangle.
\end{align*}
Similarly one can solve Equation (\ref{tm}) by the same technique for the eigenvalue $\lambda_2(p)$, but this time rewriting
out the left hand side and right hand side of the equation in terms of $X_2$. The reader may check this yields
\begin{align*}
-\lambda_2\Big(\frac{\cos(\beta)}{\sin(\beta)} +\frac{\sin(\beta)}{\cos(\beta)}\Big) = &
\frac{\sin(\beta)}{\cos(\beta)}\langle A_{\xi}(J_1\xi),J_1\xi \rangle \\
& -\frac{\cos(\beta)}{\sin(\beta)}\langle A_{\xi}(J_1Z), J_1Z\rangle.
\end{align*}

Let $M_{t}$ denote the hypersurface which is the tube of radius $t$ around $M$ obtained by traveling along the normal
geodesic $C_{\xi}$. Then $\alpha$ is constant along $C_{\xi}$ as the ambient space is symmetric by Gray's theorem. Hence
along $C_{\xi}$ we have $ J\xi(t) = \cos(\alpha)J_1\xi(t) + \sin(\alpha)J_1Z(t)$. Obviously, if $\xi(t)$ denotes the unit
normal vector field we can solve to find $\lambda_1(t)$ and $\lambda_2(t)$. Both these equations holding simultaneously for
all $t$ is equivalent to
\begin{equation}\label{la1}
\langle A_{\xi}(J_1Z(t)), J_1Z(t)\rangle =\lambda_1(t)\Big(1 - \frac{\cos^2(\beta)}{\sin^2(\beta)}\Big) +
\frac{\cos^2(\beta)}{\sin^2(\beta)}\langle A_{\xi}(J_1\xi(t)), J_1\xi(t)\rangle
\end{equation}
\begin{equation}\label{la2}
\langle A_{\xi}(J_1Z(t)), J_1Z(t)\rangle = \lambda_2(t)\Big(1 + \frac{\sin^2(\beta)}{\cos^2(\beta)}\Big) +
\frac{\sin^2(\beta)}{\cos^2(\beta)}\langle A_{\xi}(J_1\xi(t)), J_1\xi(t)\rangle.
\end{equation}
Suppose that $\cos(\beta) \neq \sin(\beta)$. Equating $\langle A_{\xi}(J_1Z(t)), J_1Z(t)\rangle$ in these two equations and
performing a routine calculation gives
$$
\lambda_1(t)\Big(1 - \frac{\cos^2(\beta)}{\sin^2(\beta)}\Big) - \lambda_2(t)\Big(1 + \frac{\sin^2(\beta)}{\cos^2(\beta)}\Big) =
(-\cos(2\beta))\langle A_{\xi}(J_1\xi(t)), J_1\xi(t)\rangle.
$$

As $\langle A_{\xi}(J_1\xi(t)), J_1\xi(t)\rangle$ is assumed to be constant with respect to $t$ one obtains $\lambda_1(t) =
c\lambda_2(t)$ for some constant $c$.

However, as $M_t$ is a tube around $M$ we can solve for $\lambda_1(t)$ and $\lambda_2(t)$ by solving the Riccati equations
explicitly. Choose a point $q$ lying at distance $t$ from $p$ along the geodesic $C_{\xi}$. This yields
$$
\lambda_1'(t) = (\lambda_1(t))^2 -4 (1+\cos(\alpha)),
$$
$$
\lambda_2'(t) = (\lambda_2(t))^2 -4 (1-\cos(\alpha)).
$$
Since $\lambda_1'(t) = c\lambda_2'(t)$, solving the above two equations forces $\lambda_2(t)$ to be constant along $C_{\xi}$
and hence $\lambda_2'(t) = 0$. But this is impossible, as $0<\alpha < \frac{\pi}{2}$. If $\cos(\beta)= \sin(\beta)$, the above
equations implies $\lambda_2(t)=0$ for all $t$, and again the Riccati equation gives a contradiction.

Therefore our supposition must be false, and hence $K_{\xi}$ must have one of the other two eigendecompositions given in the
list in \cite{berndt2}. The reader may read off from these tables that either $J\xi \in \mathcal{J}$ or $J\xi(p) \perp
\mathcal{J}(p)$, and hence $M$ has singular normal Jacobi operator at every point. This implies that $\cos(\alpha) \in
\lbrace 0, \frac{\pi}{2}\rbrace,$ a contradiction.
\endproof

For the non-compact dual of $G_2(\mathbb{C}^{m+2})$, namely the symmetric space $SU_{2,m}/S(U_2U_m)$, the same calculation
yields

\begin{thm}
There are no curvature-adapted hypersurfaces of \linebreak $SU_{2,m}/S(U_2U_m)$ if  $\cos(\alpha(p)) \notin \lbrace 0,
\frac{3}{5}, \frac{4}{5},1 \rbrace$ and  either $\langle A_{\xi}(J_1\xi(t)), J_1\xi(t)\rangle$, or  $\langle
A_{\xi}(J_1Z(t)), J_1Z(t)\rangle$, or their ratio, is constant along the normal geodesic $C_{\xi}(t)$ through all points
$p\in M$.
\end{thm}

We remark that real hypersurfaces with $\cos(\alpha)\in \lbrace 0, \frac{\pi}{2}\rbrace$ are classified in a recent paper
\cite{berndtsuh2} up to one possible exception. They obtain the list
\begin{rlist}
\item a tube of radius $r\in \mathbb{R}^+$ around a totally geodesic \\ $SU_{2,m-1}/SU_2SU_{m-1}$, or
\item a tube of radius $r\in \mathbb{R}^+$ around a totally geodesic $\mathbb{H}H^n$, where $n=2m$, or
\item a horosphere whose centre at infinity is singular, or
\item the normal space $\nu(M)$ of $M$ consists of singular tangent vectors $X$ of the form $JX \perp \mathcal{J}X$.
\end{rlist}
The same calculation as before shows these the first two families of hypersurfaces are curvature-adapted. It is conjectured
that there are no hypersurfaces with constant principal curvatures in the fourth case. This is related to a major open
problem in the submanifold geometry of symmetric spaces, which is how to find a better understanding of the geometry of
horospheres in noncompact symmetric spaces.

\textbf{Acknowledgments} This work was undertaken as part of a Ph.D. at University College Cork supported by the Irish
Research Council for Science, Engineering and Technology. The author wishes to thank his advisor J\"{u}rgen Berndt for his
valuable advice and encouragement, as well as Robert Bryant and Benjamin McKay for several illuminating conversations on the
Cayley planes.

\end{document}